\documentclass[12pt]{article}

\usepackage{amsfonts}

\usepackage{amsmath,amssymb}

 \numberwithin{equation}{subsection}

\begin{document}

\title{Some properties on $G$-evaluation and its applications to $G$-martingale decomposition }
\author{Yongsheng Song\\
\small Academy of Mathematics and Systems Science, \\
\small Chinese Academy of Sciences, Beijing, China;\\
\small yssong@amss.ac.cn}

\date{}

\maketitle

\begin{abstract}

In this article, a sublinear expectation induced by $G$-expectation
is introduced, which is called $G$-evaluation for convenience. As an
application, we prove that for any $\xi\in L^\beta_G(\Omega_T)$ with
some $\beta>1$ the decomposition theorem holds and that any
$\beta>1$ integrable symmetric $G$-martingale can be represented as
an It$\hat{o}$ integral w.r.t $G$-Brownian motion.  As a byproduct,
we prove a regular property for $G$-martingale: Any $G$-martingale
$\{M_t\}$  has a quasi-continuous version.

\end{abstract}


{\bf MSC-classification}: 60G07, 60G20, 60G44, 60G48, 60H05

\maketitle
\section{Introduction }
Recently, [P06], [P08] introduced the notion of sublinear
expectation space, which is a generalization of probability space.
One of the most important sublinear expectation space is
$G$-expectation space. As the counterpart of Wiener space in the
linear case, the notions of $G$-Brownian motion, $G$-martingale, and
It$\hat{o}$ integral w.r.t $G$-Brownian motion were also introduced.
These notions have very rich and interesting new structures which
nontrivially generalize the classical ones.

Because of the Sublinearity, the fact of $\{M_t\}$ being a
$G$-martingale does not imply that $\{-M_t\}$ is a $G$-martingale. A
surprising fact is that there exist nontrivial processes which are
continuous, decreasing and are also $G$-martingales. [P07]
conjectured  that for any $\xi\in L^1_G(\Omega_T)$, we have the
following representation:
$$X_t :=\hat{E}_t(\xi)= \hat{E}(\xi) + \int_0^tZ_sdB_s-K_t, \ t\in[0, T],$$ with $K_0=0$ and $\{K_t\}$ an
increasing process.

[P07] proved the conjecture for cylindrical functions
$L_{ip}(\Omega_T)$(see Theorem 2.18) by It$\hat{o}'s$ formula in the
setting of $G$-expectation space. So the left question is to extend
the representation to the completion $L^1_G(\Omega_T)$ of
$L_{ip}(\Omega_T)$ under norm $\|\xi\|_{1,G}=\hat{E}(|\xi|)$.

[STZ09] make a progress in this direction. They define a much
stronger norm $\|\xi\|_{{\cal L}^0_2}=\{\hat{E}[\sup_{t\in[0,
T]}\hat{E}_t(|\xi|^2)]\}^{1/2}$ on $L_{ip}(\Omega_T)$ and
generalized the above result to the completion ${\cal L}^0_2\subset
L^2_G(\Omega_T)$. The shortcoming of this result is that no
relations between the two norms are given.  The space ${\cal L}^0_2$
is just an abstract completion and we have no idea about the set
$L^2_G(\Omega_T)\backslash {\cal L}^0_2$.

The purpose of this article is to improve  the decomposition theorem
given  in [P07] and [STZ09].  The main results of the article
consist of:

 We introduce a sublinear expectation called $G$-evaluation and
investigate its properties. By presenting an estimate, which can be
seen as the substitute of Doob's maximal inequality, we proved that
for any $\xi\in L^\beta_G(\Omega_T)$ with some $\beta>1$
 the decomposition theorem holds and that any
$\beta>1$ integrable symmetric $G$-martingale can be represented as
an It$\hat{o}$ integral w.r.t $G$-Brownian motion.

As a byproduct, we prove a regular property for $G$-martingale: Any
$G$-martingale $\{M_t\}$  has a quasi-continuous version(see
Definition 5.1). We also give several estimates for variables in the
decomposition theorem, which may be useful in the follow-up work of
$G$-martingale theory.

This article is organized as follows: In section 2, we recall some
basic notions and results of $G$-expectation and the related space
of random variables. In section 3, we introduce the notion of
$G$-evaluation  and present an estimate, which can be seen as the
substitute of Doob's maximal inequality. In section 4, we prove that
for  any $\beta>1$ integrable $G$-martingale, the decomposition
theorem holds and that any $\beta>1$ integrable symmetric
$G$-martingale can be represented as an It$\hat{o}$ integral w.r.t
$G$-Brownian motion. In section 5, we prove a regular property for
$G$-martingale.

\section{Preliminary }
We present some preliminaries in the theory of sublinear
expectations and the related $G$-Brownian motions. More details of
this section can be found in [P07].

\subsection{G-expectation }

\noindent {\bf Definition 2.1} Let $\Omega$
 be a given set and let ${\cal H}$ be a linear space of real valued
functions defined on $\Omega$
 with $c \in {\cal H}$ for all constants $c$. ${\cal H}$ is considered as the
space of  ¡°random variables¡±. A sublinear expectation $\hat{E}$ on
${\cal H}$ is a functional $\hat{E}: {\cal H}\rightarrow R $
satisfying the following properties: for all $X, Y \in {\cal H}$, we
have

(a) Monotonicity: If $X\geq Y$ then $\hat{E}(X) \geq \hat{E} (Y)$.

(b) Constant preserving: $\hat{E} (c) = c$.

(c) Sub-additivity: $\hat{E}(X)-\hat{E}(Y) \leq \hat{E}(X-Y)$.

(d) Positive homogeneity: $\hat{E} (\lambda X) = \lambda \hat{E}
(X)$, $\lambda \geq 0$.

\noindent$(\Omega, {\cal H}, \hat{E})$ is called a sublinear
expectation space.

\noindent {\bf Definition 2.2} Let $X_1$ and $X_2$ be two
$n$-dimensional random vectors defined respectively in sublinear
expectation spaces $(\Omega_1, {\cal H}_1, \hat{E}_1)$ and
$(\Omega_2, {\cal H}_2, \hat{E}_2)$. They are called identically
distributed, denoted by $X_1 \sim X_2$, if $\hat{E}_1[\varphi(X_1)]
= \hat{E}_2[\varphi(X_2)]$, $\forall \varphi\in C_{l, Lip}(R^n)$,
where $ C_{l, Lip}(R^n)$ is the space of real continuous functions
defined on $R^n$ such that $$|\varphi(x) - \varphi(y)| \leq C(1 +
|x|^k + |y|^k)|x - y|, \forall x, y \in R^n,$$ where $k$ depends
only on $\varphi$.

\noindent {\bf Definition 2.3} In a sublinear expectation space
$(\Omega, {\cal H}, \hat{E})$ a random vector $Y = (Y_1,
\cdot\cdot\cdot, Y_n)$, $Y_i \in {\cal H}$ is said to be independent
to another random vector $X = (X_1, \cdot\cdot\cdot, X_m)$, $X_i \in
{\cal H}$ under $\hat{E}(\cdot)$ if for each test function
$\varphi\in C_{l, Lip}(R^m\times R^n)$ we have $\hat{E}[\varphi(X, Y
)] = \hat{E}[\hat{E} [\varphi(x, Y )]_{x=X}]$.

\noindent {\bf Definition 2.4} ($G$-normal distribution) A
d-dimensional random vector $X = (X_1, \cdot\cdot\cdot,X_d)$ in a
sublinear expectation space $(\Omega, {\cal H}, \hat{E})$ is called
$G$-normal distributed if for each $a, b\in R$ we have $$aX +
b\hat{X}\sim \sqrt{a^2 + b^2}X,$$  where $\hat{X}$ is an independent
copy of $X$. Here the letter $G$ denotes the function $$G(A) :=
\frac{1 }{2}\hat{ E}[(AX,X)] : S_d \rightarrow R,$$  where $S_d$
denotes the collection of $d\times d$ symmetric matrices.

The function $G(\cdot) : S_d \rightarrow R$ is a monotonic,
sublinear mapping on $S_d$ and $G(A) = \frac{1 }{2}\hat{
E}[(AX,X)]\leq \frac{1 }{2}|A|\hat{ E}[|X|^2]=:\frac{1
}{2}|A|\bar{\sigma}^2$ implies that there exists a bounded, convex
and closed subset $\Gamma\subset S_d^+$ such that
$$G(A)=\frac{1 }{2}\sup_{\gamma\in \Gamma}Tr(\gamma A).$$ If there exists some $\beta>0$ such that
$G(A)-G(B)\geq \beta Tr(A-B)$ for any $A\geq B$, we call the $G$-normal
distribution is non-degenerate, which is the case we consider
throughout this article.

\noindent {\bf Definition 2.5} i) Let $\Omega_T=C_0([0, T]; R^d)$
with the supremum norm, $ {\cal H}^0_T:=\{\varphi(B_{t_1},...,
B_{t_n})| \forall n\geq1, t_1, ..., t_n \in [0, T], \forall \varphi
\in C_{l, Lip}(R^{d\times n})\}$, $G$-expectation is a sublinear
expectation defined by
$$\hat{E}[\varphi( B_{t_1}-B_{t_0}, B_{t_2}-B_{t_1} ,
\cdot\cdot\cdot, B_{t_m}- B_{t_{m-1}} )]$$$$ = \tilde{E}
[\varphi(\sqrt{t_1-t_0}\xi_1, \cdot\cdot\cdot, \sqrt{t_m
-t_{m-1}}\xi_m)],$$ for all $X=\varphi( B_{t_1}-B_{t_0},
B_{t_2}-B_{t_1} , \cdot\cdot\cdot, B_{t_m}- B_{t_{m-1}} )$, where
$\xi_1, \cdot\cdot\cdot, \xi_n$ are identically distributed
$d$-dimensional $G$-normal distributed random vectors in a sublinear
expectation space $(\tilde{\Omega}, \tilde{\cal H},\tilde{ E})$ such
that  $\xi_{i+1}$ is independent to $(\xi_1, \cdot\cdot\cdot,
\xi_i)$ for each $i = 1, \cdot\cdot\cdot,m$. $(\Omega_T, {\cal
H}^0_T, \hat{E})$ is called a $G$-expectation space.

ii) For $t\in [0, T]$ and $\xi=\varphi(B_{t_1},..., B_{t_n})\in
{\cal H}^0_T$, the conditional expectation defined by(there is no
loss of generality, we assume $t=t_i$) $$\hat{E}_{t_i}[\varphi(
B_{t_1}-B_{t_0}, B_{t_2}-B_{t_1} , \cdot\cdot\cdot, B_{t_m}-
B_{t_{m-1}} )]$$$$=\tilde{\varphi}( B_{t_1}-B_{t_0}, B_{t_2}-B_{t_1}
, \cdot\cdot\cdot, B_{t_i}- B_{t_{i-1}} ),$$ where
$$\tilde{\varphi}(x_1, \cdot\cdot\cdot, x_i)=\hat{E}[\varphi( x_1,
\cdot\cdot\cdot,x_i, B_{t_{i+1}}- B_{t_{i}}, \cdot\cdot\cdot,
B_{t_m}- B_{t_{m-1}} )].$$

Let $\|\xi\|_{p, G}=[\hat{E}(|\xi|^p)]^{1/p}$ for $\xi\in{\cal
H}^0_T$ and $p\geq1$, then $\forall t\in[0, T]$, $\hat{E}_t(\cdot)$
is a continuous mapping on ${\cal H}^0_T$ with norm $\|\cdot\|_{1,
G}$ and therefore can be extended continuously to the completion
$L^1_G(\Omega_T)$ of ${\cal H}^0_T$ under norm $\|\cdot\|_{1, G}$.

\noindent {\bf Proposition 2.6} Conditional expectation defined
above has the following properties: for $X, Y \in L^1_G(\Omega_T)$

i) If $X \geq Y$ , then $\hat{E}_t(X) \geq \hat{E}_t(Y)$.

ii) $\hat{E}_t(\eta) = \eta$, for  $\eta \in L^1_G(\Omega_t)$.

iii) $\hat{E}_t(X)- \hat{E}_t(Y) \leq \hat{E}_t(X-Y)$.

iv) $\hat{E}_t(\eta X) = \eta^+\hat{E}_t (X) + \eta^-\hat{E}_t(-X)$,
for each bounded $\eta \in L^1_G(\Omega_t)$.

v) $\hat{E}_s(\hat{E}_t (X)) = \hat{E}_{t\wedge s}(X)$, in
particular, $\hat{E}(\hat{E}_t(X)) = \hat{E} (X)$.

vi) For each $X\in L^1_G(\Omega ^t_T )$ we have $\hat{E}_t (X) =
\hat{E} (X)$.

\noindent {\bf Theorem 2.7}([DHP08]) There exists a tight subset
${\cal P}\subset {\cal M}_1(\Omega_T)$ such that
$$\hat{E}(\xi)=\max_{P\in {\cal P}}E_P(\xi) \ \ \textrm{for \
all} \ \xi\in{\cal H}^0_T.$$ ${\cal P}$ is called a set that
represents $\hat{E}$.

\noindent {\bf Remark 2.8} i) [HP09] gave a new proof to the above
theorem. From the proof we know that any sublinear expectation
${\cal E}(\cdot)$ on ${\cal H}^0_T$ satisfying $${\cal
E}[(B_t-B_s)^{2n}]\leq d_n (t-s)^n, \forall n\in N,$$ has the above
representation.

ii) Let ${\cal A}$ denotes the sets that represent $\hat{E}$. ${\cal
P}^*=\{P\in {\cal M}_1(\Omega_T)| E_P(\xi)\leq \hat{E}(\xi), \
\forall \ \xi\in{\cal H}^0_T\}$ is obviously the maximal one, which
is  convex and weak compact. However, by Choquet capacitability
Theorem, all capacities induced by weak compact sets of
probabilities in ${\cal A}$ are the same, i.e. $c_{\cal
P}:=\sup_{P\in{\cal P}}P=\sup_{P\in{\cal P}'}P=:c_{{\cal P}'}$ for
any weak compact set ${\cal P}, {\cal P}'\in{\cal A}$. So we call it
the capacity induced by $\hat{E}$. In fact, By Choquet
capacitability Theorem, it suffices to prove the compact sets case.
For any compact set $K\subset\Omega_T$, there exists an decreasing
sequence $\{\varphi_n\}\subset C_b^+(\Omega_T)$ such that
$1_K\leq\varphi_n\leq1$ and $\varphi_n\downarrow 1_K$. Then by
Theorem 28 in [DHP08],$$ c_{\cal
P}(K)=\lim_{n\rightarrow\infty}\sup_{P\in{\cal
P}}E_P(\varphi_n)=\lim_{n\rightarrow\infty}\hat{E}(\varphi_n)=\lim_{n\rightarrow\infty}\sup_{P\in{\cal
P}'}E_P(\varphi_n)= c_{{\cal P}'}(K).$$

iii) All capacities induced by  sets of probabilities in ${\cal A}$
are the same on open sets. In fact, for any ${\cal P}\in {\cal A}$,
let $\bar{{\cal P}}$ be the weak closure of ${\cal P}$. Since
$c_{\bar{{\cal P}}}=c_{\cal P}$ on open sets, we get the desired
result by ii).

iv)Let $(\Omega^0, \{{\cal F}^0_t\}, {\cal F}, P^0 )$ be a filtered
probability space, and $\{W_t\}$ be a d-dimensional Brownian motion
under $P^0$.  [DHP08] proved that $${\cal P}'_M:=\{P_0\circ X^{-1}|
X_t=\int_0^th_sdW_s, h\in L^2_{\cal F}([0,T]; \Gamma^{1/2}) \}\in
{\cal A},$$ where $\Gamma^{1/2}:=\{\gamma^{1/2}| \gamma\in \Gamma\}$
and $\Gamma$ is the set in the representation of $G(\cdot)$.

v) Let ${\cal P}_M$ be the weak closure of ${\cal P}'_M$. Then under
each $P\in {\cal P}_M$, the canonical process $B_t(\omega)=\omega_t$
for $\omega\in\Omega_T$ is a martingale. In fact, for any $P\in
{\cal P}_M$, there exists $\{P_n\}\subset{\cal P}'_M$ such that
$P_n\rightarrow P$ weakly. For any $0\leq s\leq t \leq T$ and
$\varphi\in C_b(\Omega_s)$,
$E_{P_n}[\varphi(B)B_s]=E_{P_n}[\varphi(B)B_t]$ since $\{B_t\}$  is
a martingale under $P_n$. Then by the integrability of $B_t, B_s$
and the weak convergence of  $\{P_n\}$ we have
$$E_{P}[\varphi(B)B_s]=\lim_{n}E_{P_n}[\varphi(B)B_s]=\lim_{n}E_{P_n}[\varphi(B)B_t]=E_{P}[\varphi(B)B_t].$$
Thus we get the desired result. $\Box$

\noindent {\bf Definition 2.9} i) Let $c$ be the capacity induced by
$\hat{E}$. A map $X$ on $\Omega_T$ with values in a topological
space is said to be quasi-continuous w.r.t $c$ if
$$\forall \varepsilon>0, \ there \ exists \ an \ open \ set \ O \ with \ c(O)<\varepsilon \ such \ that \ X|_{O^c} \ is \ continuous.$$

ii) We say that $X: \Omega_T\rightarrow R$ has a quasi-continuous
version if there exists a quasi-continuous function $Y:
\Omega_T\rightarrow R$ with $X=Y$, $c$-q.s.. $\Box$

By the definition of quasi-continuity and iii) in Remark 2.8, we
know that the collections of quasi-continuous functions w.r.t.
capacities induced by any set(not necessary weak compact) that
represents $\hat{E}$ are the same.

 Let $\|\varphi\|_{p,G}=[\hat{E}(|\varphi|^p)]^{1/p}$ for $\varphi\in C_b(\Omega_T)$, the
 completions of
$C_b(\Omega_T)$, ${\cal H}^0_T$ and $L_{ip}(\Omega_T)$ under $\|\cdot\|_{p,G}$ are the same and denoted by $L^p_G(\Omega_T)$,
where $$L_{ip}(\Omega_T):=\{\varphi(B_{t_1},..., B_{t_n})| \forall
n\geq1, t_1, ..., t_n \in [0, T], \forall \varphi \in C_{b,
Lip}(R^{d\times n})\}$$ and $C_{b, Lip}(R^{d\times n})$ denotes the
set of bounded Lipschitz functions on $R^{d\times n}$.

\noindent {\bf Theorem 2.10}[DHP08] For  $p\geq1$ the completion
$L^p_G(\Omega_T)$ of $C_b(\Omega_T)$ is $$L^p_G(\Omega_T)=\{X\in
L^0: X \ has \ a \ q.c. \ version, \
\lim_{n\rightarrow\infty}\hat{E}[|X|^p1_{\{|X|>n\}}]=0\},$$ where
$L^0$ denotes the space of all R-valued measurable functions on
$\Omega_T$.

\subsection{Basic notions on stochastic calculus
in sublinear expectation space}

Now we shall introduce some basic notions on stochastic calculus in
sublinear expectation space $(\Omega_T, L^1_G, \hat{E})$. The
canonical process $B_t(\omega)=\omega_t$ for $\omega\in\Omega_T$ is
called $G$-Brownian motion.

\textit{For convenience of description, we only give the definition
of It$\hat{o}$ integral with respect to 1-dimensional $G$-Brownian
motion. However, all results in sections 3-5 of this article hold
for the $d$-dimensional case.}

For $p\geq1$, let $M^{p,0}_G(0, T)$ be the collection of processes
in the following form: for a given partition $\{t_0,
\cdot\cdot\cdot, t_N\} = \pi_T$ of $[0, T]$,  $$ \eta_t(\omega) =
\sum^{N-1}_{j=0} \xi_j(\omega)1_{[t_j ,t_{j+1})}(t),$$ where
$\xi_i\in L^p_G(\Omega_ {t_i})$, $i = 0, 1, 2, \cdot\cdot\cdot,
N-1$. For each $\eta\in M^{p,0}_G(0, T)$, let
$\|\eta\|_{M^{p}_G}=\{\hat{E}\int_0^T|\eta_s|^pds\}^{1/p}$ and
denote $M^{p}_G(0, T)$ the completion of $M^{p,0}_G(0, T)$ under
norm $\|\cdot\|_{M^{p}_G}$.

\noindent {\bf Definition 2.11} For each $\eta\in M^{2,0}_G(0, T)$
with the form $$\eta_t(\omega) = \sum^{N-1}_{j=0}
\xi_j(\omega)1_{[t_j,t_{j+1})}(t),$$ we define $$I(\eta) =\int_0^T
\eta(s)dB_s := \sum^{N-1}_{j=0} \xi_j(B_{t_{j+1} }-B_{t_j} ).$$

The mapping $I: M^{2,0}_G(0, T)\rightarrow L^2_G(\Omega_T)$ is
continuous and thus can be continuously extended to $M^2_G(0, T)$.

We denote for some $0\leq s \leq t \leq T$, $\int_s^t\eta_udB_u :=
\int^T_0 1_{[s,t]}(u)\eta_udB_u$. We have the following properties:

Let $\eta, \theta\in M^2_G(0, T)$ and let $0 \leq s \leq r\leq t\leq
T$. Then in $L^1_G(\Omega_T )$ we have

(i) $\int^t_s \eta_udB_u = \int^r_s \eta_udB_u+ \int^t_r
\eta_udB_u$,

(ii)$\int^t_s (\alpha\eta_u+\theta_u)dB_u=\alpha\int^t_s
\eta_udB_u+\int^t_s \theta_udB_u$, if $\alpha$ is bounded and in
$L^1_G(\Omega_s)$,

(iii)$\hat{E}_t(X +\int_t^T\eta_sdB_s) = \hat{E} (X)$, $\forall X\in
L^1_G(\Omega^t_T)$ and $\eta\in M^{2}_G(0, T)$.

\noindent {\bf Definition 2.12} Quadratic variation process of
$G$-Brownian motion defined by
$$\langle B\rangle_t = B^2_t-2\int_0^t B_sdB_s $$ is a continuous,
nondecreasing process.

For $\eta\in M^{1,0}_G(0, T)$, define $Q_{0,T} (\eta) = \int^T_0
\eta(s)d \langle B\rangle_s := \sum^{N-1}_{j=0} \xi_j(\langle
B\rangle _{t_{j+1}}- \langle B\rangle _{t_j} ) : M^{1,0}_G (0, T)
\rightarrow L^1_G(\Omega_T).$ The mapping is continuous and can be
extended to $M^{1}_G (0, T)$ uniquely.

\noindent {\bf Definition 2.13} A process $\{M_t\}$ with values in
$L^1_G(\Omega_T)$ is called a $G$-martingale if $\hat{E}_s(M_t)=M_s$
for any $s\leq t$. If $\{M_t\}$ and  $\{-M_t\}$ are both
$G$-martingale, we call $\{M_t\}$ a symmetric $G$-martingale.

\noindent {\bf Definition 2.14}  For two process $\{X_t\}, \{Y_t\}$
with values in $L^1_G(\Omega_T)$, we say $\{X_t\}$ is a version of $
\{Y_t\}$ if
$$X_t=Y_t, \ \ q.s. \ \ \forall t\in[0, T].$$

 By the same arguments as in the classical linear case, for which we
refer to [HWY92] for instance, we have the following lemma.

\noindent {\bf Lemma 2.15} Any symmetric $G$-martingale
$\{M_t\}_{t\in[0, T]}$  has a RCLL(right continuous with left limit)
version. $\Box$

In the rest of this article, we only consider the RCLL versions of
symmetric $G$-martingales.

\noindent {\bf Theorem 2.16 [P07]} For each $x \in R$, $Z\in M^2_G
(0, T)$ and  $\eta\in M^1_G (0, T)$, the process $$M_t = x +
\int_0^tZ_sdB_s + \int_0^t\eta_sd\langle B\rangle_s -\int_0^t
2G(\eta_s)ds, \ t\in[0, T]$$ is a martingale. $\Box$

\noindent {\bf Remark 2.17} Specially, $-K_t=\int_0^t\eta_sd\langle
B\rangle_s -\int_0^t 2G(\eta_s)ds$ is a $G$-martingale, which is a
surprising result because $-K_t$ is a continuous, non-increasing
process. [P07] conjectured that any $G$-martingale has the above
form and gave the following result. $\Box$

\noindent {\bf Theorem 2.18 [P07]} For all
$\xi=\varphi(B_{t_1}-B_{t0}, \cdot\cdot\cdot,
B_{t_n}-B_{t_{n-1}})\in L_{ip}(\Omega_T)$, we have the following
representation: \begin {eqnarray}\xi = \hat{E}(\xi) +
\int_0^TZ_tdB_t + \int_0^T\eta_td\langle B\rangle_t -\int_0^T
2G(\eta_t)dt.
\end {eqnarray} where $Z\in M^2_G
(0, T)$ and  $\eta\in M^1_G (0, T)$.

[STZ09] defined $\|\xi\|_{{\cal L}^0_2}=\{\hat{E}[\sup_{t\in[0,
T]}\hat{E}_t(|\xi|^2)]\}^{1/2}$ on $L_{ip}(\Omega_T)$ and
generalized the above result to the completion ${\cal L}^0_2$ of
$L_{ip}(\Omega_T)$ under $\|\cdot\|_{{\cal L}^0_2}$.

\noindent {\bf Theorem 2.19 [STZ09]} For all $\varphi\in {\cal
L}^0_2$, there exists $\{Z_t\}_{t\in[0, T]}\in M^2_G(0, T)$ and a
continuous increasing process $\{K_t\}_{t\in[0, T]}$ with $K_0=0,
K_T\in L^2_G(\Omega_T)$ and $\{-K_t\}_{t\in[0, T]}$ a $G$-martingale
such that
\begin {eqnarray}
X_t:=\hat{E}_t(\varphi)=\hat{E}(\varphi)+\int_0^tZ_sdB_s-K_t=:M_t-K_t,
\ q.s.
\end {eqnarray}
$\Box$

\section{$G$-evaluation }
In this section, we introduce an sublinear expectation which is
induced by $G$-expectation and investigate some of its properties.

For $\xi\in {\cal H}^0_T$, let ${\cal E}(\xi)=\hat{E}[\sup_{u\in[0,
T]}\hat{E}_u(\xi)]$ for all $\xi \in {\cal H}^0_T$ where  $\hat{E}$
is the $G$-expectation. For convenience, we call ${\cal E}$
$G$-evaluation. First we give the following representation for
$G$-evaluation, which is similar to that of $G$-expectation.

\noindent {\bf Theorem 3.1} There exists a weak compact subset
${\cal P}^{\cal E}\subset {\cal M}_1(\Omega)$ such that $${\cal
E}(\xi)=\max_{P\in {\cal P}^{\cal E}}E_P(\xi) \ \ \textrm{for \ all}
\ \xi\in{\cal H}^0_T.$$

{\bf Proof.}
1. Obviously, $(\Omega, {\cal H}^0_T, {\cal E})$ is a
sublinear expectation space. Then there exists a family of positive
linear functionals $ {\cal I}$ on ${\cal H}^0_T$ such that
$${\cal
E}(\xi)=\max_{I\in {\cal I}}I(\xi) \ \ \textrm{for \ all} \
\xi\in{\cal H}^0_T.$$

2. In the following, we give some calculations.

For any $0\leq s\leq t\leq T$ and $u\in[0, T]$,

$${\hat{E}_u|B_t-B_s|^{2n}}\leq\left\{\begin{array}{ll}
{|B_t-B_s|^{2n}},& \ \ \textrm{if} \ \ u\geq t,\\
{c_n(t-s)^n},& \ \ \textrm{if} \ \ u\leq
s,\\{2^{2^n-1}[c_n(t-u)^n+|B_u-B_s|^{2n}]},& \ \ \textrm{if} \ \
s<u< t.\end{array} \right.$$ Thus $\hat{E}_u|B_t-B_s|^{2n}\leq
2^{2^n-1}[c_n(t-s)^n+\sup_{u\in [s, t]}|B_u-B_s|^{2n}]$, and by
B-D-G inequality

${\cal E}|B_t-B_s|^{2n}\leq
2^{2^n-1}[c_n(t-s)^n+b_n(t-s)^n]=:d_n(t-s)^n$.

3. Noting the discussion in Remark 2.8,  we can prove the desired
representation by just the same arguments as in [HP09]. $\Box$

For $p\geq 1$ and $ \xi\in {\cal H}^0_T$, define $\|\xi\|_{p, {\cal
E}}=[{\cal E}(|\xi|^p)]^{1/p}$ and denote $L^p_{{\cal E}}(\Omega_T)$
the completion of $ {\cal H}^0_T$ under $\|\cdot\|_{p, {\cal E}}$.

We shall give an  estimate between the two norms $\|\cdot\|_{p,
{\cal E}}$ and $\|\cdot\|_{p, G}$. As the substitute of Doob's
maximal inequality, the estimate will play a critical role in the
proof to the martingale decomposition theorem in the next section.
First, we shall give a lemma.

\textit{For convenience, we say $\xi$ is symmetric if $\xi\in L^1_G$
with $\hat{E}(\xi)+\hat{E}(-\xi)=0$.}

\noindent {\bf Lemma 3.2} For $\xi\in L_{ip}(\Omega_T)$, there
exists nonnegative $K_T\in L^1_G(\Omega_T)$ such that $\xi+K_T$ is
 symmetric. Moreover, for any $1<\gamma<\beta$, $\gamma\leq2$,
$K_T\in L^\gamma_G(\Omega_T)$ and
$$\|K_T\|^\gamma_{L^\gamma_G}\leq14C^\gamma_{\beta/\gamma}
\|\xi\|^\beta_{L^\beta_G},$$  where
$C_{\beta/\gamma}=\sum_{i=1}^{\infty}i^{-\beta/\gamma}$.

{\bf Proof.} Let $\xi^n=(\xi\wedge n)\vee(-n)$ and
$\eta^n=\xi^{n+1}-\xi^n$ for $n\geq0$. Then by Theorem 2.18, for
each n, we have the following representation (2.2.1):
$$X^n_t:=\hat{E}_t(\eta^n)=M^n_t-K^n_t,$$ where $\{M^n_t\}$ is a
symmetric $G$-martingale with $M^n_T\in L^2_G(\Omega_T)$ and
$\{K^n_t\}$ is a continuous increasing process  with $K^n_0=0,
K^n_T\in L^2_G(\Omega_T)$ and $\{-K^n_t\}$ a $G$-martingale.  Fix
$P\in{\cal P}_M$. By It$\hat{o}'s$ formula
$$(\eta^n)^2=2\int_0^TX^n_tdX^n_t+[M^n]_T, \ P-a.s.$$ Take
expectation under $P$, we have
$$
E_P[(M^n_T)^2] \leq E_P[(\eta^n)^2]+2E_P(K^n_T).
$$
Take supremum over ${\cal P}_M$, we have
$$
\hat{E}[(M^n_T)^2] \leq\hat{E}[(\eta^n)^2]+2\hat{E}(K^n_T) \leq
5\hat{E}(|\eta^n|).
$$
Therefore,
$$\hat{E}[(K_T^n)^2]\leq2(\hat{E}[(\eta^n)^2]+\hat{E}[(M^n_T)^2])\leq
12\hat{E}(|\eta^n|).$$ Consequently, for any $1<\gamma<\beta$ and
$\gamma\leq2$
$$\hat{E}[(K_T^n)^\gamma]\leq\hat{E}(K_T^n)+\hat{E}[(K_T^n)^2])\leq
14\hat{E}(|\eta^n|).$$
\begin {eqnarray*}& &\hat{E}[(\sum_{i=n+1}^{n+m}K^i_T)^\gamma]\\
&\leq & [\sum_{i=n+1}^{n+m}i^{-\beta/\gamma}]^{\gamma-1}
\sum_{i=n+1}^{n+m}i^{\beta/\gamma^*}
\hat{E} [(K^i_T)^\gamma]\\
&\leq&
14C^{\gamma-1}_{\beta/\gamma}(n,m)\sum_{i=n+1}^{n+m}i^{\beta/\gamma^*}
\hat{E} (|\eta^i|)\\
&\leq& 14C^{\gamma-1}_{\beta/\gamma}(n,m)\sum_{i=n+1}^{n+m}i^{\beta/\gamma^*} c(|\xi|>i)\\
&\leq& 14\hat{E}(|\xi|^\beta)C^{\gamma}_{\beta/\gamma}(n,m).
\end {eqnarray*}
where $C_{\beta/\gamma}(n,m)=\sum_{i=n+1}^{n+m}i^{-\beta/\gamma},
\gamma^*=\gamma/(\gamma-1).$

So $\{\sum_{n=0}^NK^n_T\}$ is a Cauchy sequence in
$L^\gamma_G(\Omega_T)$. Let $K_T:=\lim_{L^\gamma_G,
N\rightarrow\infty}\sum_{n=0}^NK^n_T$, then
$\|K_T\|^\gamma_{L^\gamma_G}\leq14C^\gamma_{\beta/\gamma}
\|\xi\|^\beta_{L^\beta_G}$. Since $\eta^n+K_T^n$ is  symmetric
for each $n\geq0$, $\xi^N+\sum_{n=0}^{N-1}K^n_T$ is  symmetric
for each $N\geq1$. Consequently, $\xi+K_T$ is  symmetric.
$\Box$

\noindent {\bf Theorem 3.3} For any $\alpha\geq1$ and $\delta>0$,
$L^{\alpha+\delta}_G(\Omega_T)\subset L^\alpha_{\cal E}(\Omega_T)$.
More precisely, for any $1<\gamma<\beta:=(\alpha+\delta)/\alpha$,
$\gamma\leq2$, we have
$$\|\xi\|_{\alpha, {\cal E}}^\alpha\leq \gamma^*\{\|\xi\|^\alpha_{\alpha+\delta,
G}+14^{1/\gamma}C_{\beta/\gamma}\|\xi\|^{(\alpha+\delta)/\gamma}_{\alpha+\delta,
G}\}, \ \forall \xi\in L_{ip}(\Omega_T),$$  where
$C_{\beta/\gamma}=\sum_{i=1}^{\infty}i^{-\beta/\gamma},
\gamma^*=\gamma/(\gamma-1).$

{\bf Proof.} For $\xi\in L_{ip}(\Omega_T)$, $|\xi|^\alpha\in
L_{ip}(\Omega_T)$. By Lemma 3.2, there exists $K_T\in
L^1_G(\Omega_T)$ such that for any $1<\gamma<\beta$, $\gamma\leq2$
$K_T\in L^\gamma_G(\Omega_T)$ and $M_T:=|\xi|^\alpha+K_T$ is
symmetric. Let $M_t=\hat{E}_t(M_T)$. Then
\begin {eqnarray*}\|\xi\|_{\alpha, {\cal E}}^\alpha&=&\hat{E}[\sup_{t\in[0, T]}\hat{E}_t(|\xi|^\alpha)]\\
&\leq& \hat{E}(\sup_{t\in[0, T]}M_t)\\
&\leq&  [\hat{E}(\sup_{t\in[0, T]}M^\gamma_t)]^{1/\gamma}\\
&\leq& \gamma^*\|M_T\|_{\gamma, G}\\
&\leq&
\gamma^*\{\||\xi|^\alpha\|_{\gamma, G}+\|K_T\|_{\gamma, G}\}\\
&\leq& \gamma^*\{\|\xi\|^\alpha_{\alpha+\delta,
G}+14^{1/\gamma}C_{\beta/\gamma}\|\xi\|^{(\alpha+\delta)/\gamma}_{\alpha+\delta,
G}\}.
\end {eqnarray*} $\Box$

Let  ${\cal P}_{\cal E}$ be weak compact subsets of ${\cal
M}_1(\Omega_T)$ which represent  ${\cal E}$. Define capacity
 $c_{\cal E}(A)=\sup_{P\in{\cal
P}_{\cal E}}(A)$. We all $c_{\cal E}$ the capacity induced by ${\cal
E}$.

By the above estimate, we can get the following equivalence between
the Choquet capacities induced by $\hat{E}$ and ${\cal E}$.

\noindent {\bf Corollary 3.4} There exists $C>0$ such that for any  set $A\in{\cal
B}(\Omega_T)$, $c(A)^2\leq c_{\cal E}(A)^2\leq C c(A).$

{\bf Proof.} By Choquet  capacitability Theorem, it suffices to
prove the compact sets case. For any compact set $K\subset\Omega_T$,
there exists an decreasing sequence $\{\varphi_n\}\subset
C_b^+(\Omega_T)$ such that $1_K\leq\varphi_n\leq1$ and
$\varphi_n\downarrow 1_K$. Let $\alpha=\delta=1$ in the above
Theorem 3.3,  there exists $1<\gamma<2$ and $C>0$, such that
$$[{\cal E}(\varphi_n)]^2\leq C \hat{E}(\varphi_n).$$

Then by Theorem 28 in [DHP08],$$ c_{\cal E}(K)^2\leq Cc(K).$$
$\Box$

\noindent {\bf Corollary 3.5}  The collections of quasi-continuous
functions on $\Omega_T$ w.r.t $c$ and $c_{\cal E}$ are the same.
$\Box$

\section{Applications to $G$-martingale decomposition }

\subsection{ Generalized It$\hat{o}$ integral }
Let $H^0_G(0, T)$ be the collection of processes in the following
form: for a given partition $\{t_0, \cdot\cdot\cdot, t_N\} = \pi_T$
of $[0, T]$,  $$ \eta_t(\omega) = \sum^{N-1}_{j=0}
\xi_j(\omega)1_{[t_j ,t_{j+1})}(t),$$ where $\xi_i\in L_{ip}(\Omega_
{t_i})$, $i = 0, 1, 2, \cdot\cdot\cdot, N-1$. For each $\eta\in
H^0_G(0, T)$ and $p\geq1$, let
$\|\eta\|_{H^{p}_G}=\{\hat{E}(\int_0^T|\eta_s|^2ds)^{p/2}\}^{1/p}$
and denote $H^{p}_G(0, T)$ the completion of $H^0_G(0, T)$ under
norm $\|\cdot\|_{H^{p}_G}$. It's easy to prove that $H^{2}_G(0,
T)=M^{2}_G(0, T)$.

\noindent {\bf Definition 4.1} For each $\eta\in H^0_G(0, T)$ with
the form $$\eta_t(\omega) = \sum^{N-1}_{j=0}
\xi_j(\omega)1_{[t_j,t_{j+1})}(t),$$ we define $$I(\eta) =\int_0^T
\eta(s)dB_s := \sum^{N-1}_{j=0} \xi_j(B_{t_{j+1} }-B_{t_j} ).$$

By B-D-G inequality, the mapping $I: H^0_G(0, T)\rightarrow
L^p_G(\Omega_T)$ is continuous under $\|\cdot\|_{H^{p}_G}$ and thus
can be continuously extended to $H^p_G(0, T)$.

\subsection{ $G$-martingale decomposition}

Let ${\cal B}_t=\sigma\{B_s| s\leq t\}$, ${\cal F}_t=\cap_{r>t}
{\cal B}_r$ and $\mathbb{F}=\{{\cal F}_t\}_{t\in[0, T]}$. $\tau:
\Omega_T\rightarrow[0, T]$ is called a $\mathbb{F}$ stopping time if
$[\tau\leq t]\in {\cal F}_t$, $\forall t\in[0, T]$.

In order to prove the more general $G$-martingale decomposition
decomposition theorem, we first introduce a famous lemma, for which
we refer to [RY94].

\noindent {\bf Definition 4.3.} A positive, adapted right-continuous
process $X$ is dominated by an increasing process $A$ with
$A_0\geq0$ if $$E[X_\tau]\leq E[A_\tau]$$ for any bounded stopping
time $\tau$.

\noindent {\bf Lemma 4.4.} If $X$ is dominated by  $A$ and  $A$ is
continuous, for any $k\in(0, 1)$
$$E[(X^*_T)^k]\leq\frac{2-k}{1-k}E[A^k_T],$$ where $X^*_T=\sup_{t\in[0,
T]}X_t$.

\noindent {\bf Theorem 4.5.}  For $\xi\in L^\beta_G(\Omega_T)$ with
some $\beta>1$, $X_t=\hat{E}_t(\xi)$, $ t\in[0, T]$ has the
following decomposition:
\begin {eqnarray*}
X_t=X_0+\int_0^tZ_sdB_s-K_t, \ q.s.
\end {eqnarray*}
 where $\{Z_t\}\in H^1_G(0, T)$  and $\{K_t\}$ is a continuous
 increasing process with $K_0=0$ and $\{-K_t\}_{t\in[0, T]}$ a
$G$-martingale. Furthermore, the above decomposition is unique and
$\{Z_t\}\in H^\alpha_G(0, T)$, $K_T\in L^\alpha_G(\Omega_T)$ for any
$1\leq\alpha<\beta$.

{\bf Proof.} For $\xi\in L^\beta_G(\Omega_T)$, there exists a
sequence $\{\xi^n\}\subset L_{ip}(\Omega_T)$ such that
$\|\xi^n-\xi\|_{\beta, G}\rightarrow0$. By Theorem 2.18, we have the
following decomposition
\begin {eqnarray*}
X^{n}_t:=\hat{E}_t(\xi^n)=X_0+\int_0^tZ^n_sdB_s-K_t:=M^n_t-K^n_t, \
q.s.
\end {eqnarray*}
where $\{Z^n_t\}\in H^2_G(0, T)$ and $\{K^n_t\}$ is a continuous
increasing process with $K^n_0=0$ and $\{-K^n_t\}_{t\in[0, T]}$ a
$G$-martingale.

Fix $P\in {\cal P}_M$, by It$\hat{o}'s$ formula,
$$(X^{n}_\tau)^2=2\int_0^\tau X^n_sdX^n_s+[M^n]_\tau, \ \forall \textmd{\ stopping \ time} \ \tau.$$
Take expectation under $P$, we have
\begin {eqnarray*}
E_P[(M^n_\tau)^2]&=&E_P[(X^{n}_\tau)^2]+2E_P(\int_0^\tau
X^n_sdK^n_s)\\
&\leq&E_P[(X^{n*}_\tau)^2]+2E_P(\int_0^\tau X^{n*}_sdK^n_s),
\end {eqnarray*} where $X^{n*}_t=\sup_{0<s\leq t}|X^n_s|$.

In the following, $C_\alpha$ will always designate a universal
constant, which may vary from line to line.

$\beta\leq2$ case.

Consequently, for any $1<\alpha<\beta$, by Lemma 4.4
\begin {eqnarray*}
E_P[(M^{n*}_T)^\alpha]&\leq& C_\alpha
\{E_P[(X^{n*}_T)^\alpha]+E_P[(X^{n*}_T)^{\alpha/2}(K^n_T)^{\alpha/2}]\}\\
&\leq& C_\alpha
\{E_P[(X^{n*}_T)^\alpha]+\{E_P[(X^{n*}_T)^{\alpha}]\}^{1/2}\{E_P[(K^n_T)^{\alpha}]\}^{1/2}\},
\end {eqnarray*} where $M^{n*}_T=\sup_{0<s\leq T}|M^n_s|$.

On the other hand, $$(K^n_\tau)^2\leq2[(X^n_\tau)^2+(M^n_\tau)^2], \
\forall \tau.$$ So
\begin {eqnarray*}
E_P[(K^n_\tau)^2]&\leq&2E_P[(X^n_\tau)^2+(M^n_\tau)^2]\\
&\leq& 2E_P[(X^{n*}_\tau)^2+(M^{n*}_\tau)^2].
\end {eqnarray*}
By this, we have $$E_P[(K^n_T)^\alpha]\leq C_\alpha
E_P[(X^{n*}_T)^\alpha+(M^{n*}_T)^\alpha].$$ So
\begin {eqnarray*}
& &E_P[(M^{n*}_T)^\alpha]\\
&\leq& C_\alpha
E_P[(X^{n*}_T)^\alpha]+C_\alpha\{E_P[(X^{n*}_T)^{\alpha}]\}^{1/2}\{E_P[(M^{*n}_T)^{\alpha}]\}^{1/2}\\
&\leq& 1/2C_\alpha E_P[(X^{n*}_T)^\alpha]+1/2E_P[(M^{n*}_T)^\alpha].
\end {eqnarray*}

 Now, we have $E_P(|M^{n*}_T|^\alpha) \leq C_\alpha
 E_P[(X^{n*}_T)^\alpha]$ and $E_P(|K^n_T|^\alpha) \leq C_\alpha
 E_P[(X^{n*}_T)^\alpha]$.

Let $\widehat{X}_t:=X^n_t-X^m_t$, $\widehat{M}_t:=M^n_t-M^m_t$,
$\widehat{K}_t:=K^n_t-K^m_t$ and $\widetilde{K}_t:=K^n_t+K^m_t$. By
It$\hat{o}'s$ formula,
$$\widehat{X}^2_\tau=2\int_0^\tau\widehat{X}_sd\widehat{X}_s+[\widehat{M}]_\tau, \ \forall \tau.$$
Take expectation under $P$, we have
\begin {eqnarray*}
E_P[(\widehat{M}_\tau)^2]&=&E_P[(\widehat{X}_\tau)^2]+2E_P(\int_0^\tau
\widehat{X}_sd\widehat{K}_s)\\
&\leq&E_P[(\widehat{X}^*_\tau)^2]+2E_P(\int_0^\tau
\widehat{X}^*_sd\widetilde{K}_s),
\end {eqnarray*} where $\widehat{X}^*_t=\sup_{0<s\leq t}|\widehat{X}_s|$.

By the same arguments as above,
\begin {eqnarray*}
E_P[(\widehat{M}_T^*)^{\alpha}]&\leq& C_\alpha
\{E_P[(\widehat{X}^*_T)^\alpha]+E_P[(\widehat{X}^*_T)^{\alpha/2}(\widetilde{K}_T)^{\alpha/2}]\}\\
&\leq& C_\alpha
\{E_P[(\widehat{X}^*_T)^\alpha]+\{E_P[(\widehat{X}^*_T)^{\alpha}]\}^{1/2}\{E_P[(\widetilde{K}_T)^{\alpha}]\}^{1/2}\},
\end {eqnarray*}  where $\widehat{M}^*_T=\sup_{0<s\leq T}|\widehat{M}_s|$.

Take supremum over ${\cal P}_M$, we get $$\hat{E}(|K^n_T|^\alpha)
\leq C_\alpha \hat{E}[(X^{n*}_T)^\alpha]$$ and
$$\hat{E}[(\widehat{M}^*_T)^\alpha]\leq C_\alpha
\{\hat{E}[(\widehat{X}^*_T)^\alpha]+\{\hat{E}[(\widehat{X}^*_T)^{\alpha}]\}^{1/2}\{\hat{E}[(\widetilde{K}_T)^{\alpha}]\}^{1/2}\}.$$
By Theorem 3.3, $\sup_n\hat{E}[(X^{n*}_T)^\alpha]<\infty$ and
$\sup_{n,m \geq N}\hat{E}[(\widehat{X}^*_T)^\alpha]\rightarrow0$ as
$N$ goes to infinity. Therefore, $\sup_{n,m \geq
N}\hat{E}[(\widehat{M}^*_T)^\alpha]\rightarrow0$ and consequently
$$\sup_{n,m \geq N}\hat{E}(\sup_{t\in[0,
T]}|K^n_t-K^m_t|^\alpha)\rightarrow0$$ as $N$ goes to infinity. Then
there exists  symmetric $G$-martingale $\{M_t\}$ and a process
$\{K_t\}$ valued in $L^\alpha_G(\Omega_T)$ such
that$$\hat{E}(\sup_{t\in[0, T]}|M^n_t-M_t|^\alpha)\rightarrow0$$ and
$$\hat{E}(\sup_{t\in[0, T]}|K^n_t-K_t|^\alpha)\rightarrow0$$ as $n$ goes to infinity.

So by B-D-G inequality, there exists $\{Z_t\}\in H^\alpha_G(0, T)$
such that $\|Z-Z^n\|_{H^\alpha_G}\rightarrow0$.

Consequently,
$$X_t=\lim_{L^\alpha_G, n\rightarrow\infty}X^n_t=
\lim_{L^\alpha_G,
n\rightarrow\infty}\int_0^tZ^n_sdB_s-\lim_{L^\alpha_G,
n\rightarrow\infty}K^n_t=\int_0^tZ_sdB_s-K_t.$$

$\beta>2$ case.

 For $2<\alpha<\beta$, $$[M^n]^{\alpha/2}_T\leq C_\alpha(|X^n_T|^\alpha+|\int_0^TX^n_sdK^n_s|^{\alpha/2}+|\int_0^TX^n_sdM^n_s|^{\alpha/2}).$$
 So \begin {eqnarray*}& & E_P([M^n]^{\alpha/2}_T)\\
 &\leq&
 C_\alpha[E_P(|X^n_T|^\alpha)+E_P(|\int_0^TX^n_sdK^n_s|^{\alpha/2})+E_P(|\int_0^TX^n_sdM^n_s|^{\alpha/2})]\\
&\leq&C_\alpha\{E_P(|X^n_T|^\alpha)+\{E_P[(X^{n*}_T)^\alpha]\}^{1/2}
\{E_P[(K^{n}_T)^\alpha]\}^{1/2}\\ & & +
\{E_P[(X^{n*}_T)^\alpha]\}^{1/2}
\{E_P([M^n]^{\alpha/2}_T)\}^{1/2}\}.
 \end {eqnarray*}
On the other hand
\begin {eqnarray*} E_P[(K_T^n)^\alpha]&\leq& C_\alpha
[E_P(|X_T^n|^\alpha)+E_P(|M_T^n|^\alpha)]\\
&\leq&C_\alpha\{E_P[(X_T^{n*})^\alpha]+E_P([M^n]^{\alpha/2})\}.
\end {eqnarray*}
Therefore,  $$E_P([M^n]^{\alpha/2}_T)\leq C_\alpha
E_P[(X_T^{n*})^\alpha]$$ and $$E_P[(K_T^n)^\alpha]\leq C_\alpha
E_P[(X_T^{n*})^\alpha].$$ By the same arguments, we get

\begin {eqnarray*}
E_P([\widehat{M}]_T^{\alpha/2})\leq C_\alpha
\{E_P[(\widehat{X}^*_T)^\alpha]+\{E_P[(\widehat{X}^*_T)^{\alpha}]\}^{1/2}\{E_P[(\widetilde{K}_T)^{\alpha}]\}^{1/2}\}.
\end {eqnarray*} The rest of the proof is just similar to the
$\beta\leq2$ case.

$\Box$

\noindent {\bf Theorem 4.6} Let $\xi\in L^\beta_G(\Omega_T)$ for
some $\beta>1$ with $\hat{E}(\xi)+\hat{E}(-\xi)=0$, then there
exists $\{Z_t\}_{t\in[0, T]}\in H^1_G(0, T)$ such that
$$\xi=\hat{E}(\xi)+\int_0^TZ_sdB_s.$$ Furthermore, the above representation is unique and
$\{Z_t\}\in H^\alpha_G(0, T)$ for any $1\leq\alpha<\beta$.

{\bf Proof.} By Theorem 4.5, for $\xi\in L^\beta_G(\Omega_T)$ with
some $\beta>1$, $X_t=\hat{E}_t(\xi)$, $ t\in[0, T]$ has the
following decomposition:
\begin {eqnarray*}
\xi=\hat{E}(\xi)+\int_0^TZ_sdB_s-K_T, \ q.s.
\end {eqnarray*}
 where $\{Z_t\}\in H^\alpha_G(0, T)$ for any $1\leq\alpha<\beta$  and $\{K_t\}$ is a continuous
 increasing process with $K_0=0$ and $\{-K_t\}_{t\in[0, T]}$ a
$G$-martingale. If $\xi$ is symmetric in addition, then
$$\hat{E}(K_T)=\hat{E}(\xi)+\hat{E}(-\xi)=0.$$ So $K_T=0, q.s.$ and
$$\xi=\hat{E}(\xi)+\int_0^TZ_sdB_s.$$
$\Box$

\noindent {\bf Remark 4.7} Since  $\xi\in L^\beta_G(\Omega_T)$, we
have by B-D-G and Doob's maximal inequality
$$\{\hat{E}(\int_0^T|Z_s|^2ds)^{\beta/2}\}^{1/\beta}<\infty.$$ But
we still can't say that $\{Z_t\}_{t\in[0, T]}\in H^\beta_G(0, T)$.
Fortunately, by a stopping time technique, Corollary 5.2 in [S10]
shows that $\{Z_t\}_{t\in[0, T]}$ does belong to $H^\beta_G(0, T)$
for $\xi\in L^\beta_G(\Omega_T)$. Here we will give a direct proof
for $\beta=2$ case.

\noindent {\bf Theorem 4.8} Let $\xi\in L^2_G(\Omega_T)$ with
$\hat{E}(\xi)+\hat{E}(-\xi)=0$, then there exists $\{Z_t\}_{t\in[0,
T]}\in M^2_G(0, T)$ such that
$$\xi=\hat{E}(\xi)+\int_0^TZ_sdB_s.$$

{\bf Proof.} Let $\xi^n=(\xi\wedge n)\vee(-n)$, then
$\hat{E}[(\xi-\xi^n)^2]\rightarrow0$. Let $M^n_t=\hat{E}_t(\xi^n)$
and $-\widetilde{M}^n_t=\hat{E}_t(-\xi^n)$. By Theorem 4.5, there
exist $\{Z^n_t\}, \{\widetilde{Z}^n_t\}\in M^2_G(0, T)$  and
continuous increasing processes $\{K^n_t\}_{t\in[0, T]},
\{\widetilde{K}^n_t\}_{t\in[0, T]}$ with $K^n_0=\tilde{K}^n_0=0$ and
$\{-K^n_t\}_{t\in[0, T]}$, $\{-\widetilde{K}^n_t\}_{t\in[0, T]}$
$G$-martingales such that
\begin {eqnarray*}
M^n_t&=&M^n_0+\int_0^tZ^n_sdB_s-K^n_t=:N^n_t-K^n_t,\\
\widetilde{M}^n_t&=&\widetilde{M}^n_0-\int_0^tZ^n_sdB_s+\widetilde{K}^n_t=:\widetilde{N}^n_t+\widetilde{K}^n_t.
\end {eqnarray*}
Let
$$\hat{M}^n_t:=M^n_t-\widetilde{M}^n_t=N^n_t-\widetilde{N}^n_t-(\widetilde{K}^n_t+K^n_t)=:\hat{N}^n_t-\hat{K}^n_t.$$

Fix $P\in{\cal P}_M$. By It$\hat{o}'s$ formula,
$$0=(\hat{M}^n_T)^2=2\int_0^T\hat{M}^n_sd\hat{M}^n_s+[\hat{N}^n]_T, \ P-a.s. \ \forall t\in[0, T].$$
Take expectation in the above equation,
\begin {eqnarray*}& &E_P[(\hat{N}^n_T)^2]=
2E_P(\int_0^T\hat{M}^n_sd\hat{K}^n_s)\leq 4n E_P[\hat{K}^n_T].
\end {eqnarray*}
So $\hat{E}[(\hat{N}^n_T)^2]\leq 4n \hat{E}[\hat{K}^n_T]$. Noting
that
\begin {eqnarray*}& &\hat{E}[\hat{K}^n_T]\\
&=&\hat{E}[\hat{N}^n_T]+\hat{E}[-\hat{M}^n_T]\\
&=&\hat{E}[\xi^n]+\hat{E}[-\xi^n]\\
&=&\hat{E}[\xi^n-\xi]+\hat{E}[\xi]+\hat{E}[-(\xi^n-\xi)]-\hat{E}[\xi]\\
&\leq&2\hat{E}[|\xi^n-\xi|]\\
&\leq&2\hat{E}[|\xi|1_{[|\xi|>n]}],
\end {eqnarray*}
we get $$\hat{E}[(\hat{N}^n_T)^2]\leq 4n
\hat{E}[\hat{K}^n_T]\leq8n\hat{E}[|\xi|1_{[|\xi|>n]}]\leq8\hat{E}[|\xi|^21_{[|\xi|>n]}]\rightarrow0.$$
So
$$\hat{E}[(K^n_T)^2]\leq\hat{E}[(\hat{K}^n_T)^2]=\hat{E}[(\hat{N}^n_T)^2]\rightarrow0.$$
Let $X^n:=\xi^n-\xi=N^n_T-\xi-K^n_T$. Then
$$\hat{E}[(N^n_T-\xi)^2]\leq2\{\hat{E}[(X^n)^2]+\hat{E}[(K^n_T)^2]\}=2\{\hat{E}[(\xi^n-\xi)^2]+\hat{E}[(K^n_T)^2]\}\rightarrow0.$$

 Since $\{\eta \in L^2_G(\Omega_T)| \
\eta=\int_0^TZ_sdB_s \ \textmd{for \ some} \ Z\in M^2_G(0, T)\}$ is
closed in $L^2_G(\Omega_T)$, we proved the desired result. $\Box$

By the proof in Theorem 4.5, we can get the following estimates,
which may be useful in the follow-up work of $G$-martingale  theory.

\noindent {\bf Corollary 4.9} For $\xi, \xi'\in L^\beta_G(\Omega_T)$
with some $\beta>1$, let $\xi=M_T-K_T$ and $\xi'=M_T'-K_T'$ be the
decomposition in Theorem 4.5. Then for any $1<\alpha<\beta$,
$1<\gamma<\beta/\alpha, \gamma\leq2$, there exists $C_{\alpha,
\beta, \gamma}$ such that
$$\|K_T\|_{\alpha, G}^\alpha \leq C_{\alpha, \beta, \gamma}
\{\|\xi\|^\alpha_{\beta, G}+\|\xi\|^{\beta/\gamma}_{\beta, G}\}$$
and \begin {eqnarray*}& &\|M_T-M'_T\|_{\alpha, G}^\alpha \leq
C_{\alpha, \beta, \gamma} \{\|\xi-\xi'\|^\alpha_{\beta,
G}+\|\xi-\xi'\|^{\beta/\gamma}_{\beta,
G}\}\\
&+&C_{\alpha, \beta, \gamma} \{\|\xi-\xi'\|^\alpha_{\beta,
G}+\|\xi-\xi'\|^{\beta/\gamma}_{\beta,
G}\}^{1/2} \{1+\|\xi\|^{\beta/\gamma}_{\beta,
G}+\|\xi'\|^{\beta/\gamma}_{\beta,
G}\}^{1/2}.
\end {eqnarray*}

\section{Regular properties for $G$-martingale}

\noindent {\bf Definition 5.1}  We say that a process $\{M_t\}$ with
values in $L^1_G(\Omega_T)$ is quasi-continuous if

\emph{$\forall \varepsilon>0$, there exists open set $G$ with
$c(G)<\varepsilon$ such that $M_\cdot(\cdot)$ is continuous on
$G^c\times[0, T]$}.

\noindent {\bf Corollary 5.2} Any $G$-martingale $\{M_t\}$ with
$M_T\in L^\beta_G(\Omega_T)$ for some $\beta>1$ has a
quasi-continuous version.

{\bf Proof.}  For $\xi\in L_{ip}(\Omega_T)$, $M_t=\hat{E}_t(\xi)$ is
continuous on $[0, T]\times\Omega_T$. In fact, for
$\xi=\varphi(B_{t_1}-B_{t_0}, \cdot\cdot\cdot, B_{t_n}-B_{t_n-1})$
with $\varphi\in C_{b, lip}(R^n)$, $M_t(\cdot)$ is obvious
continuous on $\Omega_T$ for fixed $t\in[0, T]$. On the other hand,
fix $\omega\in \Omega_T$
\begin {eqnarray*}
& &|M_{t_n}(\omega)-M_{t_{n-1}}(\omega)|\\
&\leq& |\varphi(x_1, \cdot\cdot\cdot,x_{n-1},
x_n)-\hat{E}[\varphi(x_1, \cdot\cdot\cdot,x_{n-1},
B_{t_n}-B_{t_n-1})]|\\
&\leq& L \hat{E}(|B_{t_n}-B_{t_n-1}-x_n|)\\
&\leq& L
(\bar{\sigma}(t_n-t_{n-1})^{1/2}+|\omega_{t_n}-\omega_{t_{n-1}}|),
\end {eqnarray*} where $L$ is the Lipschitz constant of $\varphi$
and $x_i=B_{t_i}(\omega)-B_{t_{i_1}}(\omega)$. In fact, the above
estimate holds for any $s, t\in[t_{i-1}, t_i]$  for some 1$\leq
i\leq n$. Then for any $s_k\rightarrow s$ and $\omega^k\rightarrow
\omega$, \begin {eqnarray*}
& &|M_{s_k}(\omega^k)-M_s(\omega)|\\
&\leq&
|M_{s_k}(\omega^k)-M_s(\omega^k)|+|M_{s}(\omega^k)-M_s(\omega)|\\
&\leq& L (\bar{\sigma}|s_k-s|^{1/2}+|\omega^k_{s_k}-\omega^k_s|)
+|M_{s}(\omega^k)-M_s(\omega)|\rightarrow0.
\end {eqnarray*}

For $\xi=M_T\in L^\beta_G(\Omega_T)$, by Theorem 3.3, there exists
$1<\alpha<\beta$ and $\{\xi^n\}\subset L_{ip}(\Omega_T)$ such that
$\|\xi^n-\xi\|_{{\cal E}, \alpha}\rightarrow0$. Let
$M_t^n:=\hat{E}_t(\xi^n)$. Then $$\sup_{m>n}\hat{E}[\sup_{t\in[0,
T]}(|M^m_t-M^n_t|]\leq\sup_{m>n}\|\xi^n-\xi^m\|_{{\cal E},
\alpha}\rightarrow0$$ as $n$ goes to infinity. So there exists a
subsequence $\{n_k\}$ such that $$\sup_{m>n_k}\hat{E}[\sup_{t\in[0,
T]}(|M^m_t-M^{n_k}_t|]<1/4^{k}.$$ Consequently,
$$\hat{E}[\sum_{k=1}^\infty \sup_{t\in[0,
T]}(|M^{n_{k+1}}_t-M^{n_k}_t|)]\leq\sum_{k=1}^\infty\hat{E}[
\sup_{t\in[0, T]}(|M^{n_{k+1}}_t-M^{n_k}_t|)]<\infty.$$ Then there
exists $\{\widetilde{M}_t\}$ such that $$\sup_{t\in[0,
T]}(|M^{n_{k}}_t-\widetilde{M}_t|)\leq \sum_{i=k}^\infty
\sup_{t\in[0, T]}(|M^{n_{i+1}}_t-M^{n_i}_t|), \ \forall k\geq1.$$
For any $\varepsilon>0$, let $O^\varepsilon_k:=[\sup_{t\in[0,
T]}(|M^{n_{k+1}}_t-M^{n_{k}}_t|)>1/(2^k\varepsilon)]$ and
$O^\varepsilon=\cup_{k=1}^\infty O^\varepsilon_k$. Then
$c(O^\varepsilon)\leq\sum_{k=1}^\infty
c(O^\varepsilon_k)<\varepsilon$ and on $(O^\varepsilon)^c$

$$\sup_{t\in[0,
T]}(|M^{n_{k}}_t-\widetilde{M}_t|)\leq \sum_{i=k}^\infty
\sup_{t\in[0, T]}(|M^{n_{i+1}}_t-M^{n_i}_t|)\leq
1/(2^{k-1}\varepsilon), \ \forall k\geq1.$$   So
$$\sup_{\omega\in(O^\varepsilon)^c}\sup_{t\in[0, T]}(|M^{n_{k}}_t-\widetilde{M}_t|)\rightarrow0$$
and $\{\widetilde{M}_t\}$ is a quasi-continuous version of
$\{M_t\}$. $\Box$

\noindent {\bf Theorem 5.3} Any $G$-martingale $\{M_t\}$  has a
quasi-continuous version.

{\bf Proof.}   Let $\xi:=M_T$ and $\xi^n=(\xi\wedge n)\vee(-n)$. For
$m>n$, let $\{X^n_t\}, \{X^m_t\}, \{X^{n,m}_t\}$ be the
quasi-continuous versions of $\{\hat{E}_t(\xi^n)\},
\{\hat{E}_t(\xi^m)\}$, $\{\hat{E}_t(|\xi^n-\xi^m|)\}$ respectively.

We claim that
$$\hat{E}[\sup_{m>n}\sup_{t\in[0,
T]}(|X^m_t-X^n_t|\wedge1)]\leq\hat{E}[\sup_{m>n}\sup_{t\in[0,
T]}(X^{n,m}_t\wedge1)]\downarrow0.$$ Otherwise, there exists
$\varepsilon>0$ such that $\hat{E}[\sup_{m>n}\sup_{t\in[0,
T]}(X^{n,m}_t\wedge1)]>\varepsilon$ for all $n\in N$. Consequently,
for each $n\in N$, there exists $m(n)>n$ such that
$\hat{E}[\sup_{t\in[0, T]}(X^{n,m(n)}_t\wedge1)]>\varepsilon$.

 Noting that
$$\varepsilon<\hat{E}[\sup_{t\in[0, T]}(X^{n,m(n)}_t\wedge1)]\leq
c(\sup_{t\in[0, T]}X^{n,m(n)}_t>\varepsilon/2)+\varepsilon/2,$$ we
have $c(\sup_{t\in[0, T]}X^{n,m(n)}_t>\varepsilon/2)>\varepsilon/2$.
Since
$$[\sup_{t\in[0, T]}X^{n,m(n)}_t>\varepsilon/2]=\pi (\{(\omega, t)| \
X^{n,m(n)}_t(\omega)>\varepsilon/2\}),$$ the projection of
$\{(\omega, t)| \ X^{n,m(n)}_t(\omega)>\varepsilon/2\}$ on $\Omega$,
we have stopping time $\tau_n\leq T$ such that $\hat{E}(
X^{n,m(n)}_{\tau_n})>\varepsilon^2/4$ by section theorem.

On the other hand, by Theorem 4.5, $\{X^{n,m(n)}_t\}$ has the
following decomposition $$X^{n,m(n)}_t=M^n_t-K^n_t,$$ where
$\{M^n_t\}$ is a symmetric $G$-martingale and $\{-K^n_t\}$ a
negative $G$-martingale with $M^n_T, K^n_T\in L^1_G(\Omega_T)$. By
Theorem 4.5 and Corollary 5.2, $\{M^n_t\}, \{-K^n_t\}$ can be taken
to be quasi-continuous.

So we have
 as $n\rightarrow\infty$
\begin {eqnarray*}
& &\hat{E}(X^{n,m(n)}_{\tau_n})\\
&\leq& \hat{E}(M^n_{\tau_n})+\hat{E}(-K^n_{\tau_n})\\
&\leq& \hat{E}(M^n_T)\\
&\leq&
\hat{E}(X^{n,m(n)}_T)+\hat{E}(K^n_T)\\
&\leq&2\hat{E}(X^{n,m(n)}_T) \rightarrow0.
\end {eqnarray*}
This is a contradiction.

Therefore, by the same arguments as in Corollary 5.2, $\{M_t\}$  has
a quasi-continuous version. $\Box$




\providecommand{\bysame}{\leavevmode\hbox
to3em{\hrulefill}\thinspace}
\providecommand{\MR}{\relax\ifhmode\unskip\space\fi MR }
\providecommand{\MRhref}[2]{%
  \href{http://www.ams.org/mathscinet-getitem?mr=#1}{#2}
} \providecommand{\href}[2]{#2}

\end{document}